\documentclass[a4paper,12pt]{amsart}
\usepackage{amssymb}
\usepackage{ifthen}
 \usepackage[dvips]{graphicx}
\nonstopmode \numberwithin{equation}{section}
\usepackage{amssymb}
\usepackage{ifthen}
\usepackage{graphicx}
\usepackage{amsmath}
\usepackage[T1]{fontenc} 

\nonstopmode \numberwithin{equation}{section}
\theoremstyle{plain}
\newtheorem{thm}[equation]{Theorem}
\newtheorem{cor}[equation]{Corollary}
\newtheorem{lem}[equation]{Lemma}
\newtheorem{prop}{Proposition}

\newtheorem{conj}{Conjecture}
\newtheorem{innercustomthm}{Theorem}
\newenvironment{customthm}[1]
  {\renewcommand\theinnercustomthm{#1}\innercustomthm}
  {\endinnercustomthm}

\theoremstyle{definition}

\newtheorem{prob}{Problem}
\newtheorem{rem}{Remark}[section]

\newcounter{minutes}
\divide\time by 60
\newcounter{hours}
\multiply\time by 60
\addtocounter{minutes}{-\time}

\newcounter {own}
\def\theown {\thesection       .\arabic{own}}

\newenvironment{pf}[1][]{%
 \vskip 3mm
 \noindent
 \ifthenelse{\equal{#1}{}}%
  {{\slshape Proof. }}%
  {{\slshape #1.} }%
 }%
{\qed\bigskip}

\newcounter{alphabet}
\newcounter{tmp}

\def\be{\begin{equation}}
\def\ee{\end{equation}}

\newcommand{\bee}{\begin{enumerate}}
\newcommand{\eee}{\end{enumerate}}

\newcommand{\blem}{\begin{lem}}
\newcommand{\elem}{\end{lem}}
\newcommand{\bthm}{\begin{thm}}
\newcommand{\ethm}{\end{thm}}
\newcommand{\bcor}{\begin{cor}}
\newcommand{\ecor}{\end{cor}}
\newcommand{\beg}{\begin{examp}}
\newcommand{\eeg}{\end{examp}}
\newcommand{\begs}{\begin{examples}}
\newcommand{\eegs}{\end{examples}}
\newcommand{\bdefe}{\begin{defin}}
\newcommand{\edefe}{\end{defin}}
\newcommand{\bprob}{\begin{prob}}
\newcommand{\eprob}{\end{prob}}
\newcommand{\bei}{\begin{itemize}}
\newcommand{\eei}{\end{itemize}}

\newcommand{\bcon}{\begin{conj}}
\newcommand{\econ}{\end{conj}}
\newcommand{\bcons}{\begin{conjs}}
\newcommand{\econs}{\end{conjs}}
\newcommand{\bprop}{\begin{prop}}
\newcommand{\eprop}{\end{prop}}
\newcommand{\br}{\begin{rem}}
\newcommand{\er}{\end{rem}}
\newcommand{\brs}{\begin{rems}}
\newcommand{\ers}{\end{rems}}
\newcommand{\bo}{\begin{obser}}
\newcommand{\eo}{\end{obser}}
\newcommand{\bos}{\begin{obsers}}
\newcommand{\eos}{\end{obsers}}
\newcommand{\bpf}{\begin{pf}}
\newcommand{\epf}{\end{pf}}
\newcommand{\ba}{\begin{array}}
\newcommand{\ea}{\end{array}}
\newcommand{\beq}{\begin{eqnarray}}
\newcommand{\beqq}{\begin{eqnarray*}}
\newcommand{\eeq}{\end{eqnarray}}
\newcommand{\eeqq}{\end{eqnarray*}}

\begin{document}

\title{On Some Subclass of Harmonic Close-to-convex  Mappings}

\author{Nirupam Ghosh}
\address{Nirupam Ghosh,
Department of Mathematics,
Indian Institute of Technology Kharagpur,
Kharagpur-721 302, West Bengal, India.}
\email{nirupam91@iitkgp.ac.in}

\author{A. Vasudevarao}
\address{A. Vasudevarao,
Department of Mathematics,
Indian Institute of Technology Kharagpur,
Kharagpur-721 302, West Bengal, India.}
\email{alluvasu@maths.iitkgp.ernet.in}

\subjclass[{AMS} Subject Classification:]{Primary 30C45, 30C80}

\def\thefootnote{}
\footnotetext{ {\tiny File:~\jobname.tex,
printed: \number\year-\number\month-\number\day,
          \thehours.\ifnum\theminutes<10{0}\fi\theminutes }
} \makeatletter\def\thefootnote{\@arabic\c@footnote}\makeatother

\begin{abstract}
Let $\mathcal{H}$ denote the class of harmonic functions $f$ in $\mathbb{D}:= \{z\in \mathbb{C}:|z| < 1\}$ normalized by $f(0) = 0 = f_z(0) -1$. For $\alpha \geq 0$, we consider  the following class
$$\mathcal{W}^0_{\mathcal{H}}(\alpha):= \{f = h + \overline{g}\in\mathcal{H}: {\rm Re\,}(h'(z) + \alpha z h''(z)) >|g'(z) + \alpha z g''(z)|, \quad z\in \mathbb{D}\}.
$$
In this paper, we first prove the coefficient conjecture of Clunie and Sheil-Small for functions in the class $\mathcal{W}^0_{\mathcal{H}}(\alpha)$.
We also prove  growth theorem, convolution, convex combination properties for functions in  the class $\mathcal{W}^0_{\mathcal{H}}(\alpha)$.
Finally, we  determine the value of $r$ so that the partial sums of  functions in the class $\mathcal{W}^0_{\mathcal{H}}(\alpha)$ are close-to-convex in $|z|<r$.
\end{abstract}

\maketitle
\pagestyle{myheadings}
\markboth{ Nirupam Ghosh and A. Vasudevarao }{On Some Subclass of Harmonic Close-to-convex  Mappings}

\section{Introduction}

Let $\mathcal{H}$ be the class of complex-valued harmonic functions $f$ in the unit disk $\mathbb{D} = \{z\in \mathbb{C}:|z| < 1\}$ normalized by $f(0) = 0 = f_z(0) - 1$. Any function $f$ in $\mathcal{H}$ has the canonical representation $f = h + \overline{g}$, where
\begin{equation}\label{niru-p2-e001}
h(z) = z + \sum_{n = 2}^{\infty}a_n z^n \quad \mbox{and} \quad g(z) = \sum_{n = 1}^{\infty}b_n z^n.
\end{equation}
Here both $h$ and $g$ are  analytic functions in  $\mathbb{D}$ and are called analytic and co-analytic part of $f$ respectively. In particular, for $g(z) = 0$, the class $\mathcal{H}$  reduces to the class $\mathcal{A}$,  consisting of analytic functions in $\mathbb{D}$ with $f(0)= 0$ and $f'(0) = 1$. If $f = h + \overline{g}$ then the Jacobian $J_f(z)$ of $f$ is defined by $ J_f(z) = |h'(z)|^2 - |g'(z)|^2$ and we say $f$ is sense preserving  if $J_f(z) > 0$ in $\mathbb{D}$. Let $\mathcal{S_H}$ be the subclass of $\mathcal{H}$ consisting of univalent (that is, one-to-one) and sense preserving harmonic mappings. If $g(z) = 0$ in $\mathbb{D}$  then the class $\mathcal{S_H}$  reduces to the class $\mathcal{S}$, containing univalent analytic functions in $\mathbb{D}$ with $f(0)= 0$ and $f'(0) = 1$. Furthermore, if $f = h + \overline{g} \in \mathcal{S_H}$ then $|g'(0)| =|b_1|<1$ (as $J_f(0) = 1 - |g'(0)|^2 =1 - |b_1|^2 >0$). Thus  function
$$
F(z) = \frac{f - b_1\overline{f}}{1 - |b_1|^2}
$$
belongs to the class $\mathcal{S_H}$. Clearly $F$ is univalent  because $F$ is an affine mapping of $f$.
A simple observation shows that $F_{\overline{z}}(0) = 0$.
Thus we may  restrict our attention to the following subclass
$$
\mathcal{S}^0_{\mathcal{H}} := \{ f\in \mathcal{S_H} : b_1 = \overline{f_{\overline{z}}(0)} = 0 \}.
$$
Hence for any function $f = h + \overline{g} $ in $\mathcal{S}^0_{\mathcal{H}}$, its analytic and co-analytic parts can be represented by
\begin{equation}\label{niru-p2-e001a}
h(z) = z + \sum_{n = 2}^{\infty}a_n z^n \quad \mbox{and} \quad g(z) = \sum_{n = 2}^{\infty}b_n z^n.
\end{equation}
The family $\mathcal{S}^0_{\mathcal{H}}$ is known to be compact and normal, whereas $\mathcal{S_H}$ is normal but not compact. In 1984, Clunie and Sheil-Small \cite{Clunie-1984} 
investigated the class $\mathcal{S_H}$ and its geometric subclasses. Since then, the class $\mathcal{S_H}$ and its subclasses have been extensively studied
(see \cite{Bshouty-Lyzzaik-2011}, \cite{Bshouty-Joshi-2013}, \cite{Clunie-1984}, \cite{Kalaj-Ponnusamy-Matti-2014}, \cite{Wang-Liang-2001}).
A domain $\Omega$ is called starlike with respect to a point $z_0\in \Omega$ if the line segment joining $z_0$ to any point in
$\Omega$ lies in $\Omega$. In particular if $z_0=0$ then $\Omega$ is called starlike domain.
A complex-valued harmonic mapping $f\in\mathcal{H}$ is said to be starlike if $f(\mathbb{D})$ is starlike domain with respect to the origin. The class of harmonic starlike functions in $\mathbb{D}$ is denoted by $\mathcal{S^*_H}$. A domain $\Omega$ is said to be convex domain if it is starlike with respect to
every point in $\Omega$. A function  $f$ in  $\mathcal{H}$ is said to be convex if $f(\mathbb{D})$ is convex. The class of harmonic convex  mappings in $\mathbb{D}$ is denoted by $\mathcal{K_H}$. A domain $\Omega$ is said to be close-to-convex if the complement of $\Omega$ can be written as  union of non-intersecting half lines. A function $f\in \mathcal{H}$ is called close-to-convex in $\mathbb{D}$ if  $f(\mathbb{D})$ is close-to-convex domain. The class of  harmonic close-to-convex mappings in $\mathbb{D}$ is denoted by $\mathcal{C_H}$.
If $f_{\overline{z}}(0) = 0$ then the classes $\mathcal{S^*_H}$, $\mathcal{K_H}$ and $\mathcal{C_H}$ reduce to
$\mathcal{S^*_H}^0$, $\mathcal{K}^0_{\mathcal{H}}$ and $\mathcal{C}^0_{\mathcal{H}}$ respectively.
In analytic case, $\mathcal{S}^*$($\mathcal{K}$ and $\mathcal{C}$ respectively) is a subclasses of $\mathcal{S}$ that contains functions $f$  such that  $f(\mathbb{D})$ is starlike (convex and close-to-convex respectively) domain. 

Clinie and Sheil-Small \cite{Clunie-1984} proved the following result which gives a sufficient condition for a harmonic function $f$ to be close-to-convex.

\begin{lem}{\cite{Clunie-1984}}\label{niru-p2-l001}
Suppose $h$ and $g$ are analytic in $\mathbb{D}$ with $|g'(0)|< |h'(0)|$ and $h + \epsilon g$ is close-to-convex analytic function for each $\epsilon~ (|\epsilon| = 1)$. Then $f = h + \overline{g}$ is harmonic close-to-convex in $\mathbb{D}$.
\end{lem}


A classical problem for functions in the class $\mathcal{S}$ where functions are of the form
 $$
 f(z)  = z + \sum_{n =0}^{\infty}a_n z^n
 $$
  is to find the sharp upper bound for the absolute value of the coefficients $a_n$ for $n\geq 2$. In 1916, Bieberbach \cite{Bieberbach-1916} proved that if $f\in\mathcal{S}$ then $|a_2|\leq2$ and conjectured that $|a_n|\leq n$ for  $n\ge 2$. In 1985, de Branges \cite{Branges-1985} proved this conjecture affirmatively. The following analogous type of conjecture was proposed by Clunie and Sheil-Small \cite{Clunie-1984}  for  functions in  the class $\mathcal{S}^0_{\mathcal{H}}$.
\begin{conj}\label{niru-p2-c001}
Let $f = h + \overline{g}$ belong to $\mathcal{S}^0_{\mathcal{H}}$, where the representation of $h$ and $g$ are given by (\ref{niru-p2-e001a}). Then $|a_n|\leq A_n$ and $|b_n|\leq B_n$ for  $n\geq 2$ where
\begin{equation}\label{niru-p2-e005}
A_n = \frac{(2n + 1)(n + 1)}{6} \quad \mbox{and} \quad B_n = \frac{(2n - 1)(n - 1)}{6}.
\end{equation}
\end{conj}

Clunie and Sheil-Small \cite{Clunie-1984} verified Conjecture \ref{niru-p2-c001} for typically real functions. In 1990, Sheil-Small \cite{Sheil-Small-1990} proved Conjecture \ref{niru-p2-c001} for functions $f$ in $\mathcal{S}^0_{\mathcal{H}}$ such that  $f(\mathbb{D})$ is starlike with respect to the origin or  $f(\mathbb{D})$ is convex in one direction.
In 2001, Wang and Liang \cite{Wang-Liang-2001} proved Conjecture \ref{niru-p2-c001} for  functions in the class $\mathcal{C}^0_{\mathcal{H}}$. However Conjecture \ref{niru-p2-c001} is still open for the class $\mathcal{S}^0_{\mathcal{H}}$. Equality occurs in (\ref{niru-p2-e005}) for the harmonic Koebe function
\begin{equation}\label{niru-p2-e006}
K(z) = \frac{z - \frac{1}{2}z^2 + \frac{1}{6}z^3}{(1 - z)^3} + \overline{\frac{\frac{1}{2}z^2 + \frac{1}{6}z^3}{(1 - z)^3}}.
\end{equation}
A harmonic function $f\in\mathcal{H}$ is said to be fully starlike (resp. fully convex ) if each $|z|< r$ is mapped onto a starlike (resp. convex) domain.
In 2013, Nagpal and  Ravichandran \cite{Nagpal-Ravichandran-2013} proved the following theorem.
\begin{customthm}{A}\label{niru-p2-t001}\cite{Nagpal-Ravichandran-2013}
A sense preserving harmonic function $f = h + \overline{g}$ is fully starlike in $\mathbb{D}$ if the analytic functions $h + \epsilon g$ are starlike in $\mathbb{D}$ for each $\epsilon ~ (|\epsilon|= 1)$.

\end{customthm}

In 2013, Li and Ponnusamy \cite{Li-Ponnusamy-2013} investigated  the properties of functions in the class
$$
\mathcal{P}^0_{\mathcal{H}} := \{f = h + \overline{g}\in\mathcal{H} : {\rm Re}(h'(z))> |g'(z)|,\quad z\in \mathbb{D}\}.
$$
The class $\mathcal{P}^0_{\mathcal{H}}$ is closely related to the class $\mathcal {R} := \{f\in \mathcal{S}: {\rm Re}(f'(z))> 0, \quad z\in\mathbb{D}\}$ which was introduced by MacGregor \cite{MacGregor-1962}. It has been proved that a harmonic function $f = h + \overline{g}$ belongs to the class $\mathcal{P}^0_{\mathcal{H}}$ if and only if the analytic functions $h + \epsilon g$ belong to $\mathcal {R}$ for each $\epsilon$ ($|\epsilon| = 1$) (see \cite{Li-Ponnusamy-2013}, \cite{Ponnusamy-Li-2013}). Using this property, Li and Ponnusamy \cite{Li-Ponnusamy-2013} have obtained the coefficient bounds and radius of convexity for the functions in the class $\mathcal{P}^0_{\mathcal{H}}$.

In 1977, Chichra \cite{Chichra-1977} studied the subclass $\mathcal{W}(\alpha)$, consisting  of $f\in\mathcal{A}$ such that   ${\rm Re\,}(f'(z) + {\alpha}z f''(z))> 0$ for $z\in \mathbb{D}$. It has been proved that  for ${\rm Re\,}\alpha \geq 0,$ the members of $\mathcal{W}(\alpha)$ are univalent in $\mathbb{D}$ and $\mathcal{W}(\alpha)$ is a subclass of the class of analytic close-to-convex  functions. In 2010, the regions of variability for functions in the class $\mathcal{W}(\alpha)$ was studied by Ponnusamy and Vasudevarao \cite{Ponnusamy-Vasu-2010}.  In 1982,  R. Singh and S. Singh \cite{Singh-Singh-1981} proved that $\mathcal{W}(1)$ is a subclass of analytic starlike functions.
In 2014, Nagpal and Ravichandran \cite{Nagpal-Ravichandran-2014} studied the following class
 $$
 \mathcal{W}^0_{\mathcal{H}} := \{ f = h + \overline{g}\in\mathcal{H}: {\rm Re}(h'(z) +  z h''(z)) >|g'(z) +  z g''(z)|,\quad z\in \mathbb{D}\}
 $$
 which is the harmonic analogue of the class $\mathcal{W}(1)$.  It is known that $\mathcal{W}^0_{\mathcal{H}}$ is a subclass of $\mathcal{S^*_H}^0$ and $\mathcal{P}^0_{\mathcal{H}}$. In particular, the members of $\mathcal{W}^0_{\mathcal{H}}$ are fully starlike in $\mathbb{D}$. The sharp coefficient bounds and the growth theorem for functions in the class $\mathcal{W}^0_{\mathcal{H}}$ have been investigated in \cite{Nagpal-Ravichandran-2014}. It has been  proved that the class $\mathcal{W}^0_{\mathcal{H}}$ is closed under convolution and  convex combinations.

 .

For two analytic functions
$$
F_1(z) =  \sum _{n = 0}^{\infty} a_n z^{n} \quad \mbox{and} \quad F_2(z) =  \sum _{n = 0}^{\infty} b_n z^{n},
$$
 the convolution (or Hadamard product) is defined by
$$
F_1 * F_2 =  \sum _{n = 0}^{\infty} a_n b_n z^{n} ,\quad  z\in \mathbb{D}.
$$
Analogously, for harmonic functions $f_1 = h_1 + \overline{g_1}$ and $f_2 = h_2 + \overline{g_2}$ in $\mathcal{H}$, the convolution  is defined as
$$
f_1 * f_2 = h_1 * h_2 + \overline {g_1 * g_2}.
$$
In 1973, Ruscheweyh and Sheil-Small \cite{Ruscheweyh-Sheil-Small-1973} proved that the class of convex analytic functions is closed under convolution. 
It is well-known that if $f\in\mathcal{C}$  and $g\in\mathcal{S}^*$ then $f*g\in\mathcal{S}^*$ and the convolution of convex function with close-to-convex function is also close-to-convex. However, the class of analytic starlike functions is not closed under the convolution. In 1984, Clunie and Sheil-Small {\cite{Clunie-1984}} proved that if $f$ is harmonic convex function and $\phi$ is  analytic convex function, then the function $f * (\phi + \alpha \overline{\phi})$ is harmonic close-to-convex function for all $\alpha~ (|\alpha| < 1)$.
For the extensive study on convolution of harmonic mappings we refer to \cite{Dorff-2001}, \cite{Dorff-2012} and \cite{Goodloe-2002}.

Let  $f(z)=\sum_{n=0}^\infty a_nz^n$  be in the class $\mathcal{S}$. Then  the $n^{\rm th}$ partial sum of $f(z)$ is defined by
$$s_n(f)=\sum_{k=0}^n a_kz^k \quad\mbox { for } n\in\mathbb{N}
$$
where $a_0 = 0$ and $a_1 = 1$.
Analogously in harmonic case, the p, q-th section/partial sum of harmonic function $f = h + \overline{g}$ given by (\ref{niru-p2-e001}) is defined as follows:
$$
s_{p,q} = s_p(h) + \overline {s_q(g)}
$$
where $ s_p(h) =\sum_{k =1}^{p} a_k z^k$ and  $s_q(g) = \sum_{k = 1}^{q}{b_k z^k}$, $p,q \geq 1$ with $a_1 = 1.$

In 1928, Szeg\"o \cite{szego-1928} proved a remarkable result which asserts   that every section $s_n(f)$ of a function $f\in\mathcal{S}$ is univalent in the disk $|z|< 1/4$. The number $1/4$ is the best possible 
as is evident from the second partial sum of Koebe function $k(z) = z/ (1 - z)^2$. For $f\in \mathcal{S}$, determining the exact  radius of univalence $r_n$ of $s_n(f)$ remains an open problem. However, many  related problems concerning sections have been solved for various  geometric subclasses of $\mathcal{S}$.  In 1941, Robertson \cite{Robertson-1941} studied  the partial sums of multivalently starlike functions (see also \cite{Robertson-1936}). In 1988, Ruscheweyh \cite{Ruscheweyh-1988} proved a strong result by showing that the partial sums, $s_n(f)$  are starlike in the disk $|z|< 1/4$ not only for the functions $f$ in the class $\mathcal{S}$ but also for the closed convex hull of $\mathcal{S}$.  For many interesting results on sections of analytic functions
we refer to \cite{Obradovic-Ponnusamy-2013a}, \cite{Ponnusamy-Sahoo-Yanagihara-2014} and \cite{Silverman-1988} (also see \cite{Bharanedhar-Ponnusamy-2014, Obradovic-Ponnusamy-2014, Singh--1970}). 
In 2013, Li and Ponnusamy \cite{Ponnusamy-Li-2013} discussed the properties of sections of functions in the class
$$
 \mathcal{P}^{0}_{\mathcal{H}}(\alpha) := \{ f = h + \overline {g}: {\rm Re\,}(h'(z) - \alpha) > |g'(z)| \quad \mbox{for} \quad z \in \mathbb{D}\}.
$$
In 2015, Li and Ponnusamy \cite{Li-Ponnusamy-2015}
investigated the properties of sections of stable harmonic convex functions.

In this  paper we introduce the following class $\mathcal{W}^0_{\mathcal{H}}(\alpha)$ (for $\alpha \geq 0$)
 $$
 \mathcal{W}^0_{\mathcal{H}}(\alpha):= \{f = h + \overline{g}\in\mathcal{H}: {\rm Re\,}(h'(z) + \alpha z h''(z)) >|g'(z) + \alpha z g''(z)| \quad \mbox{for}\quad z\in \mathbb{D}\}
 $$
with  $g'(0)=0$.
The organization of this paper as follows: In section $2$, we  prove that the class $\mathcal{W}^0_{\mathcal{H}}(\alpha)$ is a subclass of close-to-convex harmonic mappings. We also obtain the sharp coefficient bounds and growth theorem for functions in the class $\mathcal{W}^0_{\mathcal{H}}(\alpha)$. In Section $3$, we obtain the convolution and convex combination properties of functions in the class $\mathcal{W}^0_{\mathcal{H}}(\alpha)$. Finally, in section $4$, we determine  $r$ so that the partial sums $s_{p, q}(f)$ of $f\in\mathcal{W}^0_{\mathcal{H}}(\alpha)$ are close-to-convex in the disk $|z|< r.$

\section{Main Result}


\begin{thm}\label{niru-p2-t005}
A harmonic mapping $f = h + \overline{g}$ is in $\mathcal{W}^0_{\mathcal{H}}(\alpha)$ if and only if the analytic function $F = h + \epsilon {g}$ belongs to $\mathcal{W}(\alpha)$ for each $|\epsilon| = 1$.
\end{thm}

\begin{pf}
If $f = h + \overline{g} \in \mathcal{W}^0_{\mathcal{H}}(\alpha)$ then for each $|\epsilon| = 1 $,
\begin{align*}
 {\rm Re}(F'(z) + \alpha z F''(z))&= {\rm Re}((h'(z) + \alpha z h''(z)) + \epsilon (g'(z) + \alpha z g''(z)))\\
 & > {\rm Re}((h'(z) + \alpha z h''(z)) - |(g'(z) + \alpha z g''(z))|\\
 & >0 \quad \mbox{for}\quad z\in\mathbb{D}.
\end{align*}
Hence $F = h + \epsilon {g} \in\mathcal{W}(\alpha)$ for each $|\epsilon| = 1$. Conversely, if $F\in\mathcal{W}(\alpha)$ then
$$
{\rm Re}((h'(z) + \alpha z h''(z) + \epsilon (g'(z) + \alpha z g''(z)) >0 \quad \mbox{for}\quad \mathbb{D}
$$
or, equivalently,
$$
{\rm Re}((h'(z) + \alpha z h''(z)) > {\rm Re}(-\epsilon(g'(z) + \alpha z g''(z))) \quad \mbox{for}\quad \mathbb{D}.
$$
Since $\epsilon~ (|\epsilon|= 1)$ is arbitrary, for an appropriate choice of $\epsilon$ we  obtain
$${\rm Re}((h'(z) + \alpha z h''(z)) > |g'(z) + \alpha z g''(z)| \quad \mbox{for} \quad z\in\mathbb{D}.
$$
Consequently  $f\in \mathcal{W}^0_{\mathcal{H}}(\alpha).$
\end{pf}

Note that for $\alpha \geq 0$, each function in the class $\mathcal{W}(\alpha)$ is close-to-convex. Hence by Lemma \ref{niru-p2-l001},  $\mathcal{W}^0_{\mathcal{H}}(\alpha)$ is a subclass of $\mathcal{C}^0_{\mathcal{H}}$ for $\alpha \geq 0$.
In 1977, Chichra \cite{Chichra-1977}  proved that  if  $0 \leq \beta< \alpha$ then $\mathcal{W}(\alpha)\subset \mathcal{W}(\beta) $. Thus $\mathcal{W}^0_{\mathcal{H}}(\alpha)\subset \mathcal{W}^0_{\mathcal{H}}(\beta) $ if $0 \leq \beta< \alpha$.
In view of this, $\mathcal{W}(\alpha)$ is starlike for $\alpha\geq 1$ because the class $\mathcal{W}(1)$ is starlike.
For $\alpha \geq 1$, if $f\in\mathcal{W}^0_{\mathcal{H}}(\alpha)$ then $h + \epsilon {g} \in\mathcal{W}(\alpha) $ is starlike function in $\mathbb{D}$ for each $\epsilon~ (|\epsilon|=1)$. Hence by Theorem \ref{niru-p2-t001}, $\mathcal{W}^0_{\mathcal{H}}(\alpha)\subset \mathcal{S^*_H}^0$ for all $\alpha \geq 1$. In particular, for $\alpha \geq 1$ each member of $\mathcal{W}^0_{\mathcal{H}}(\alpha)$ is fully starlike.

The following  theorems give sharp coefficient bounds for  functions in the class $\mathcal{W}^0_{\mathcal{H}}(\alpha)$.

\begin{thm}\label{niru-p2-t010}
Let $f = h + \overline{g} \in \mathcal{W}^0_{\mathcal{H}}(\alpha)$ for $\alpha\geq0$ be of the form (\ref{niru-p2-e001a}). Then for  $n\geq 2$,
\begin{equation}\label{niru-p2-e010}
\displaystyle|b_n|\leq \frac {1}{\alpha n^2 + n(1 - \alpha)}.
\end{equation}
The result is sharp for the function $f(z)$ which is given by $f(z) = z + \frac {1}{\alpha n^2 + n(1 - \alpha)}\overline{z^n}$.
\end{thm}

\begin{pf}
Let $f \in \mathcal{W}^0_{\mathcal{H}}(\alpha)$. Then
\begin{equation}\label{niru-p2-e015}
 {\rm Re}(h'(z) + \alpha z h''(z)) >|g'(z) + \alpha z g''(z)| \quad \mbox{for} \quad z\in\mathbb{D} .
\end{equation}
Using the series expansion of $g(z)$ and (\ref{niru-p2-e015}) we have
\begin{align*}
r^{n-1}(\alpha n^2 + (1 - \alpha) n) |b_n| &\leq \frac{1}{2 \pi}\int_{0}^{2 \pi}|g'(r e^{i \theta}) + \alpha {r e^{i \theta}} g''(r e^{i \theta})|\, d\theta\\
& \leq \frac{1}{2 \pi}\int_{0}^{2 \pi}{\rm Re}(h'(r e^{i \theta}) + \alpha (r e^{i \theta}) h''(r e^{i \theta}))\,d\theta\\
& = 1.
\end{align*}
Letting $r \rightarrow 1^{-}$ we obtain the desired bound.
To show the bound in (\ref{niru-p2-e010}) is sharp, we consider $f(z) = z + \frac {1}{\alpha n^2 + n(1 - \alpha)}\overline{z^n}$. It is easy to see that $f\in\mathcal{W}^0_{\mathcal{H}}(\alpha)$ and $|b_n(f)| = \frac {1}{\alpha n^2 + n(1 - \alpha)}.$
\end{pf}

\begin{thm}\label{niru-p2-t015}
Let $f = h + \overline{g} \in \mathcal{W}^0_{\mathcal{H}}(\alpha)$ for $\alpha\geq0$ be of the form (\ref{niru-p2-e001a}). Then for any $n\geq 2$,
\begin{enumerate}
\item[(i)] $\displaystyle |a_n| + |b_n|\leq \frac {2}{\alpha n^2 + n(1 - \alpha)}; $\\[2mm]

\item[(ii)] $\displaystyle ||a_n| - |b_n||\leq \frac {2}{\alpha n^2 + n(1 - \alpha)};$\\[2mm]

\item[(iii)] $\displaystyle |a_n|\leq \frac {2}{\alpha n^2 + n(1 - \alpha)}.$
\end{enumerate}
All these results are sharp for the function $f(z)$ which is given by $f(z) = z + \sum _{n = 2}^{\infty}\frac{2}{\alpha n^2 + n(1 - \alpha)}z^n.$
\end{thm}

\begin{pf}
It is sufficient to prove only the first inequality as the rest  follow from it.
 Since $f= h + \overline{g} \in \mathcal{W}^0_{\mathcal{H}}(\alpha)$  by Theorem \ref{niru-p2-t005}, $h + \epsilon g$ is in the class $\mathcal{W}(\alpha)$ for  $\epsilon ~(|\epsilon| =1)$. Thus for  each $|\epsilon| = 1$ we have
$$
{\rm Re\,}((h + \epsilon g)' + \alpha z (h + \epsilon g)'') > 0 \quad \mbox{for} \quad z\in\mathbb{D}.
$$
This implies there exists an analytic function $p(z)$ which is of the form $p(z) = 1 + \sum_{n = 1}^{\infty} p_n z^n$  with  ${\rm Re\,} p(z) > 0$ in $\mathbb{D}$ such that
\begin{equation}\label{niru-p2-e020}
h'(z) + \alpha z h''(z) + \epsilon (g'(z) + \alpha z g''(z)) = p(z).
\end{equation}
Comparing the coefficients on the both sides of (\ref{niru-p2-e020}) we obtain
\begin{equation}\label{niru-p2-e022}
(\alpha n^2 +(1 - \alpha)n)(a_n + \epsilon b_n) = p_{n-1}\quad \mbox{for}\quad n\geq2.
\end{equation}
Since $|p_n|\leq 2$ for  $n \geq 1$ and $\epsilon$ ($|\epsilon| = 1$) is arbitrary, it follows from (\ref{niru-p2-e022}) that
\begin{equation}\label{niru-p2-e023}
(\alpha n^2 +(1 - \alpha)n)(|a_n| + |b_n|) \leq 2
\end{equation}
which yields the first inequality.

To prove the sharpness of (\ref{niru-p2-e023}), we consider the function  $f(z) = z + \sum _{n = 2}^{\infty}\frac{2}{\alpha n^2 + n(1 - \alpha)}z^n$. It is easy to see that $f\in \mathcal{W}^0_{\mathcal{H}}(\alpha)$ and all the three inequalities are sharp.
\end{pf}

\begin{rem}
\begin{enumerate}
\item[(i)] For $\alpha = 0$, the class $\mathcal{W}^0_{\mathcal{H}}(\alpha)$ reduces to $\mathcal{R}^0_{\mathcal{H}}$ which was studied by Li and Ponnusamy \cite{Ponnusamy-Li-2013}. If we put  $\alpha = 0$ in Theorem \ref{niru-p2-t010} and Theorem \ref{niru-p2-t015}, we obtain the results
    \cite[Theorems 1 and 2]{Ponnusamy-Li-2013} as a particular case.

\item [(ii)] For $\alpha = 1$, the class $\mathcal{W}^0_{\mathcal{H}}(\alpha)$ reduces to $\mathcal{W}^0_{\mathcal{H}}(1) = \{f = h + \overline{g} : {\rm Re}(h'(z) + z h''(z))> |g'(z) + z g''(z)|\}$ which was studied by Nagpal and Ravichandran \cite{Nagpal-Ravichandran-2014}.  If we put $\alpha = 1$ in Theorem \ref{niru-p2-t010} and Theorem \ref{niru-p2-t015}, we obtain the result \cite[Theorm 3.5]{Nagpal-Ravichandran-2014} as a particular case.
\end{enumerate}
\end{rem}


\begin{thm}\label{niru-p2-t025}
Let $f = h + \overline{g} \in \mathcal{W}^0_{\mathcal{H}}(\alpha)$  be as in (\ref{niru-p2-e001a}) with  $0<\alpha \leq 1$. Then
\begin{equation}\label{niru-p2-e021}
|z| + 2\sum_{n = 2}^{\infty}\frac{(-1)^{n-1} |z|^n}{\alpha n^2 + n(1 - \alpha)} \leq |f(z)|\leq |z| + 2\sum_{n = 2}^{\infty}\frac{|z|^n}{\alpha n^2 + n(1 - \alpha)}.
\end{equation}
This result is sharp for the function $f(z) = z + \sum _{n = 2}^{\infty}\frac{2}{\alpha n^2 + n(1 - \alpha)}z^n$ and its rotations.
\end{thm}

\begin{pf}
Let $f= h + \overline{g}$ be in the class $\mathcal{W}^0_{\mathcal{H}}(\alpha)$. Then $ F = h + \epsilon g$ is in the class $\mathcal{W}(\alpha)$ and for each $|\epsilon| = 1$ we have
\begin{equation}\label{niru-p2-e023aa}
{\rm Re\,}(F'(z) + \alpha z F''(z)) > 0 \quad \mbox{for}\quad z\in\mathbb{D}.
\end{equation}
Thus there exists an analytic function $\omega(z)$  with $\omega(0)=0$ and $|\omega(z)|<1 $ in $\mathbb{D}$ such that
\begin{equation}\label{niru-p2-e024}
F'(z) + \alpha z F''(z) = \frac{1 + \omega(z)}{1 - \omega(z)}.
\end{equation}
A simplification of (\ref{niru-p2-e024}) gives
\begin{align*}\label{niru-p2-e025}
 z^{1/ \alpha} F'(z) &= \frac{1}{\alpha}\int_{0}^{z} {\xi}^{\frac{1}{\alpha} -1} \frac{1 + \omega(\xi)}{1 - \omega(\xi)}\, d\xi\\ \nonumber
& =\frac{1}{\alpha}\int_{0}^{|z|}{(t e^{i \theta})}^{\frac{1}{\alpha} -1}\frac{1 + \omega(t e^{i \theta})}{1 - \omega(t e^{i \theta})} e^{i\theta}\, dt.
\end{align*}

\noindent Therefore we have

\begin{align*}
|z^{1/ \alpha} F'(z)|&=\left|\frac{1}{\alpha}\int_{0}^{|z|}{(t e^{i \theta})}^{\frac{1}{\alpha} -1}\frac{1 +\omega(t e^{i\theta})}{1-\omega(t e^{i \theta})} e^{i\theta}\,dt \right|\\
& \leq \frac{1}{\alpha}\int_{0}^{|z|}{t }^{\frac{1}{\alpha} -1}\frac{1 + t}{1-t}\,dt
\end{align*}
and\

\begin{align*}
|z^{1/ \alpha} F'(z)|&=\left|\frac{1}{\alpha}\int_{0}^{|z|}{(t e^{i \theta})}^{\frac{1}{\alpha} -1}\frac{1 +\omega(t e^{i\theta})}{1-\omega(t e^{i \theta})} e^{i\theta}\,dt \right|\\
&\geq \frac{1}{\alpha}\int_{0}^{|z|}t^{\frac{1}{\alpha} -1}{\rm Re\,}\frac{1 +\omega(t e^{i\theta})}{1-\omega(t e^{i \theta})}\,dt\\
& \geq \frac{1}{\alpha}\int_{0}^{|z|}t^{\frac{1}{\alpha} -1}{\rm Re\,}\frac{1 - t}{1+t}\,dt.
\end{align*}
Further computation gives
\begin{equation}\label{niru-p2-e026}
|F'(z)| = |h'(z) + \epsilon g'(z)|  \leq 1 + 2 \sum_{n =1}^{\infty}\frac{|z|^n}{1 + \alpha n}.
\end{equation}
and
\begin{equation*}
|F'(z)| = |h'(z) + \epsilon g'(z)|  \geq 1 + 2 \sum_{n =1}^{\infty}\frac{(-1)^n |z|^n}{1 + \alpha n}.
\end{equation*}
Since $\epsilon~ (|\epsilon|=1)$ was arbitrary, it follows from (\ref{niru-p2-e026}) that 
\begin{equation*}\label{niru-p2-e030}
|h'(z)| +  |g'(z)| \leq 1 + 2 \sum_{n =1}^{\infty}\frac{|z|^n}{1 + \alpha n}
\end{equation*}
and
\begin{equation*}\label{niru-p2-e030}
|h'(z)| -  |g'(z)| \geq 1 + 2 \sum_{n =1}^{\infty}\frac{(-1)^n|z|^n}{1 + \alpha n}.
\end{equation*}
\noindent Let $\Gamma$ be the radial segment from $0$ to $z$, then
\begin{align*}
|f(z)| &= \left|\int_{\Gamma}\frac{\partial f}{\partial \xi}\, d\xi + \frac{\partial f}{\partial {\overline{\xi}}} \, d\overline{\xi} \right| \leq \int_{\Gamma}(|h'(\xi)| + |g'(\xi)|)|\,d\xi|\\
& \leq \int_{0}^{|z|} \left(1 + 2 \sum_{n =1}^{\infty}\frac{|t|^n}{1 + \alpha n} \right)\, dt = |z| + 2 \sum _{n = 2}^{\infty}\frac{|z|^n}{\alpha n^2 + (1 - \alpha)n}
\end{align*}
and also
\begin{align*}
|f(z)| &= \int_{\Gamma}\left|\frac{\partial f}{\partial \xi}\, d\xi + \frac{\partial f}{\partial {\overline{\xi}}} \, d\overline{\xi} \right| \geq \int_{\Gamma}(|h'(\xi)| - |g'(\xi)|)|\,d\xi|\\
& \geq \int_{0}^{|z|} \left(1 + 2 \sum_{n =1}^{\infty}\frac{(-1)^n|t|^n}{1 + \alpha n} \right)\, dt = |z| + 2 \sum _{n = 2}^{\infty}\frac{(-1)^n|z|^n}{\alpha n^2 + (1 - \alpha)n}
\end{align*}

\noindent The equality in (\ref{niru-p2-e021}) holds for the function $f(z)$ which is given by 
$$f(z) = z + \sum _{n = 2}^{\infty}\frac{2}{\alpha n^2 + n(1 - \alpha)}z^n
$$ 
and its rotations.
\end{pf}

\newpage
The following result is a sufficient condition for functions to be in the class $\mathcal{W}^0_{\mathcal{H}}(\alpha)$.
\begin{thm}
Let $f = h + \overline{g} \in \mathcal{S}^0_{\mathcal{H}}$ be of the form (\ref{niru-p2-e001a}) and satisfies the condition
\begin{equation}\label{niru-p2-e031}
\sum_{n =2}^{\infty}(\alpha n^2 + (1 - \alpha) n)(|a_n| + |b_n|)) < 1.
\end{equation}
Then $f\in\mathcal{W}^0_{\mathcal{H}}(\alpha)$.
\end{thm}

\begin{pf} Let $f=h+\overline{g}\in \mathcal{S}^0_{\mathcal{H}}$. Using the series representation of $h(z)$ given by (\ref{niru-p2-e001a}), we obtain
\begin{align}\label{niru-p2-e031aa}
{\rm Re\,}(h'(z) + \alpha z h''(z)) & = 1 + {\rm Re\,}\left(\sum_{n=2}^{\infty} (\alpha n^2 + (1 - \alpha) n) a_n z^{n-1}\right)\\\nonumber
& \geq 1 - \left|\sum_{n=2}^{\infty} (\alpha n^2 + (1 - \alpha) n) a_n z^{n-1}\right|\\\nonumber
& \geq 1 - \sum_{n=2}^{\infty}(\alpha n^2 + (1 - \alpha) n)|a_n|.
\end{align}
In view of (\ref{niru-p2-e031}), the inequality (\ref{niru-p2-e031aa})  reduces to
\begin{align*}
{\rm Re\,}(h'(z) + \alpha z h''(z)) & > \sum_{n=2}^{\infty}(\alpha n^2 + (1 - \alpha) n)|b_n|\\
& \geq \left|\sum_{n=2}^{\infty}(\alpha n^2 + (1 - \alpha) n)b_n \right|\\
& =|g'(z) + \alpha z g''(z)|
\end{align*}
and hence $f\in\mathcal{W}^0_{\mathcal{H}}(\alpha)$.
\end{pf}

\section{Convex combinations and Convolution}
In this section, we show that the class $\mathcal{W}^0_{\mathcal{H}}(\alpha)$ is closed  under convex combinations. Also we show that $\mathcal{W}^0_{\mathcal{H}}(\alpha)$ is closed under  convolution.

\begin{thm}
The class $\mathcal{W}^0_{\mathcal{H}}(\alpha)$ is closed under convex combinations.
\end{thm}

\begin{pf}
Let $f_i:=h_i + \overline{g_i}\in \mathcal{W}^0_{\mathcal{H}}(\alpha)$ for $i = 1,2,\ldots n$ and  $\sum _{i = 1}^{n} t_i = 1$ ($0\leq t_i\leq1$).
The convex combination of $f_i$'s can  be written as
$$
f(z) = \sum_{i =1}^{n} t_i f_i(z) = h(z) + \overline{g(z)}
$$
where $h(z) = \sum_{i =1}^{n} t_i h_i(z)$ and $g(z) = \sum_{i =1}^{n} t_i g_i(z)$. Then $h$ and $g$ both are analytic in $\mathbb{D}$ with $h(0) = g(0)= h'(0)-1 = g'(0) =0$. A simple computation yields
\begin{align*}
{\rm Re\,}(h'(z) + \alpha z h''(z)) &= {\rm Re\,}\left(\sum_{i = 1}^{n} t_i  (h'_i(z) + \alpha z h''_i(z))\right)\\
& > \sum_{i = 1}^{n} t_i |g'_i(z) + \alpha z g''_i(z)|\\
& \geq |g'(z) + \alpha z g''(z)|.
\end{align*}
This shows that $f\in\mathcal{W}^0_{\mathcal{H}}(\alpha)$.
\end{pf}

A sequence $\{c_n\}_{n=0}^{\infty}$ of  non-negative numbers is said to be a convex null sequence if $c_n \rightarrow {0}$ as $n \rightarrow \infty$  and
$$
c_0 - c_1 \geq c_1 - c_2 \geq c_2 - c_3 \geq ...\geq c_{n-1} - c_n \geq ...\geq 0.
$$

\noindent To prove the convolution results on the class $\mathcal{W}^0_{\mathcal{H}}(\alpha)$ we need the following lemmas:
\begin{lem}\label{niru-p2-l005}\cite{Singh-Singh-1989}
Let $\{c_n\}_{n=0}^{\infty}$ be a convex null sequence. Then the function
$$
q(z) = \frac{c_0}{2} + \sum_{n =1}^{\infty}c_n z^n
$$
is analytic and ${\rm Re\,}q(z)> 0$ in $\mathbb{D}$.
\end{lem}

\begin{lem}\label{niru-p2-l010}\cite{Singh-Singh-1989}
Let $p(z)$ be an analytic function in the unit disk $\mathbb{D}$ with $p(0) = 1 $ and ${\rm Re\,}(p(z))> 1/2$ in $\mathbb{D}$. Then for any analytic function $f$ in $\mathbb{D}$, the function $p * f$ takes values in the convex hull of the image of $\mathbb{D}$ under $f$.
\end{lem}

\noindent Using Lemma \ref{niru-p2-l005} and Lemma \ref{niru-p2-l010}, we  prove the following interesting lemma.

\begin{lem}\label{niru-p2-l015}
Let $F$ be in the class $\mathcal{W}(\alpha)$. Then ${\rm Re}(\frac{F(z)}{z})> 1/2$.
\end{lem}

\begin{pf}
 Let $F\in\mathcal{W}(\alpha)$ be given by $F(z) = z + \sum _{n = 2}^{\infty}A_n z^n.$ Since $F\in\mathcal{W}(\alpha)$, ${\rm Re\,} (F'(z) + \alpha z F''(z)) > 0$  in $\mathbb{D}$, which is equivalent to
 $$
 {\rm Re\,}\left(1 + \sum_{n=2}^{\infty}A_n(n^2 \alpha + n (1 - \alpha)) z^{n-1}\right) > 0 \quad \mbox{for} \quad z\in\mathbb{D}.
 $$
Therefore  ${\rm Re\,}p(z) > 1/2$ in $\mathbb{D}$, where
$$
p(z) = 1 + \frac{1}{2}\sum_{n=2}^{\infty}A_n(n^2 \alpha + n (1 - \alpha)) z^{n-1}.
$$
Consider a sequence $\{c_n\}_{n = 0}^{\infty}$ defined by $c_0 = 1$ and $c_{n-1} = \frac{2}{n^2\alpha + n (1 - \alpha)}$ for $n\geq 2$ . Then it is easy to see that  $c_n \rightarrow {0}$ as $n \rightarrow \infty$  and
$$
c_0 - c_1 \geq c_1 - c_2 \geq c_2 - c_3 \geq ...\geq c_{n-1} - c_n \geq ...\geq 0.
$$
Thus the sequence $\{c_n\}_{n=0}^{\infty}$ is a convex null sequence. In view of Lemma \ref{niru-p2-l005}, the function
$$
q(z) = \frac{1}{2} + \sum_{n = 2}^{\infty}\frac{2}{n^2\alpha + n (1 - \alpha)} z^{n-1}
$$
is analytic and ${\rm Re}(q(z))> 0$ in $\mathbb{D}$.
A simple computation shows that
\begin{align*}
\frac{F(z)}{z} &= 1 + \sum_{n =2}^{\infty}A_n z^{n-1}\\
& = p(z) *  \left(1 + \sum_{n = 2}^{\infty}\frac{2}{n^2\alpha + n (1 - \alpha)}z^{n-1}\right).
\end{align*}
An application of Lemma \ref{niru-p2-l010} yields ${\rm Re}\left(\frac{F(z)}{z}\right)> 1/2$ for $z\in\mathbb{D}$.
\end{pf}

\begin{lem}\label{niru-p2-l019}
Let $F_1$ and $F_2$ be in $\mathcal{W}(\alpha)$. Then the Hadamard product $F_1* F_2$ is in $\mathcal{W}(\alpha)$.
\end{lem}

\begin{pf}
Let $F_1(z) = z + \sum_{n =2}^{\infty} A_n z^n $ and $F_2(z) = z + \sum_{n =2}^{\infty} B_n z^n$.
Then the convolution of $F_1$  and $F_2$ is given by
$$
F(z) =(F_1* F_2)(z) = z + \sum_{n =2}^{\infty} A_n B_n z^n.
$$
Since
$ z F'(z) = z F_1'(z) * F_2(z),$
a computation shows that
\begin{equation}\label{niru-p2-e0035}
F'(z) + z \alpha F''(z) = (F_1'(z) + z \alpha F_1''(z) ) * \left(\frac{F_2(z)}{z} \right).
\end{equation}
Since $F_1, F_2 \in\mathcal{W}(\alpha)$, ${\rm Re\,}(F_1'(z) + z \alpha F_1''(z) ) > 0$ for $z\in\mathbb{D}$. In view of Lemma \ref{niru-p2-l015}, ${\rm Re}\left({F_2(z)}/{z} \right) > {1}/{2}$ in $\mathbb{D}$. An application of  Lemma \ref{niru-p2-l010} to (\ref{niru-p2-e0035}) gives ${\rm Re\,}(F'(z) + z \alpha F''(z) ) > 0$ in $\mathbb{D}$. Hence $F = F_1* F_2$ belongs to the class $\mathcal{W}(\alpha)$.
\end{pf}

\noindent Using Lemma \ref{niru-p2-l019} we shall prove the following interesting result.

\begin{thm}\label{niru-p2-t035}
If $f_1$ and $f_2$ are in $\mathcal{W}^0_{\mathcal{H}}(\alpha)$ then $f_1 * f_2$ is in $\mathcal{W}^0_{\mathcal{H}}(\alpha)$.
\end{thm}

\begin{pf}
Let $f_1 = h_1 + \overline{g}_1$ and $f_2 = h_2 + \overline{g}_2$ be two functions in the class $\mathcal{W}^0_{\mathcal{H}}(\alpha)$. The convolution of $f_1$ and $f_2$ is given by  $f_1 * f_2 = h_1*h_2 + \overline{g_1 * g_2}$. To show $f_1 * f_2 $ is in the class $\mathcal{W}^0_{\mathcal{H}}(\alpha)$, it is sufficient to prove that $F = h_1*h_2 + \epsilon ({g_1 * g_2})$ is in $\mathcal{W}(\alpha)$  for each $\epsilon~ (|\epsilon| = 1).$ Since by Lemma \ref{niru-p2-l019}, $\mathcal{W}(\alpha)$ is closed under convolution, for each $\epsilon~ (|\epsilon| = 1)$, $h_i + \epsilon g_i\in\mathcal{W}(\alpha)$ for $i = 1,2$. Hence the following functions
$$
F_1 = (h_1 - g_1) * (h_2 - \epsilon g_2)
$$ and
$$
F_2 = (h_1 + g_1) * (h_2 + \epsilon g_2)
$$
are in the class $\mathcal{W}(\alpha)$. Again, since $\mathcal{W}(\alpha)$ is closed under convex combination, the function
$$
\frac{1}{2}(F_1 + F_2 ) = (h_1*h_2) + \epsilon ({g_1 * g_2}) = F
$$
is in the class $\mathcal{W}(\alpha)$. Hence $\mathcal{W}^0_{\mathcal{H}}(\alpha)$ is closed under convolution.
\end{pf}

In 2002, Goodloe \cite{Goodloe-2002} considered the Hadamard product of harmonic function with an analytic function defined as follows:
$$
f \tilde{*} \phi = h* \phi + \overline{g * \phi}
$$
where $f = h + \overline{g}$ is harmonic and $\phi$ is analytic function in $\mathbb{D}.$

\begin{thm}\label{niru-p2-t040}
Let $f\in\mathcal{W}^0_{\mathcal{H}}(\alpha)$ and $\phi\in\mathcal{A}$ be such that ${\rm Re}\left(\frac{\phi(z)}{z}\right) > 1/2$  for $z\in\mathbb{D}$. Then $f \tilde{*} \phi\in \mathcal{W}^0_{\mathcal{H}}(\alpha).$
\end{thm}

\begin{pf}
Let $f = h + \overline{g}$ be in the class $\mathcal{W}^0_{\mathcal{H}}(\alpha)$. Then we have
$$
f \tilde{*} \phi = h*\phi + \overline{g * \phi}.
$$
In view of  Theorem \ref{niru-p2-t005}, to prove $f \tilde{*} \phi$ is in the class $\mathcal{W}^0_{\mathcal{H}}(\alpha)$, it suffices to prove that
$F(z) =  h*\phi + \epsilon ({g * \phi})$ belongs to the class $\mathcal{W}(\alpha)$ for each $\epsilon$ ($|\epsilon| = 1$).
Since $f\in\mathcal{W}^0_{\mathcal{H}}(\alpha)$, for each $\epsilon$ ($|\epsilon| = 1$) the function $F_1(z) = h + \epsilon g$ belongs to $\mathcal{W}(\alpha)$. Therefore  $F = F_1 * \phi$ and
$$
F'(z) + \alpha z F''(z) = (F_1'(z) + \alpha z F_1''(z)) * \frac{\phi(z)}{z}.
$$
Since ${\rm Re}\left(\frac{\phi(z)}{z}\right) > 1/2$ and  ${\rm Re}(F_1'(z) + \alpha z F_1''(z)) > 0$ in $\mathbb{D}$, in view of Lemma \ref{niru-p2-l010}, it can be seen that $F\in \mathcal{W}(\alpha)$.
\end{pf}

\noindent As a consequence of Theorem \ref{niru-p2-t040}, we obtained the following result.

\begin{cor}
Suppose $f$ belongs to $\mathcal{W}^0_{\mathcal{H}}(\alpha)$ and $\phi\in\mathcal{K}$. Then $f \tilde{*} \phi\in \mathcal{W}^0_{\mathcal{H}}(\alpha)$.
\end{cor}

\begin{pf}
Since $\phi$ is convex, ${\rm Re}\left(\frac{\phi(z)}{z}\right) > 1/2$ for $z$ in $\mathbb{D}$ and hence the desired result follows  from Theorem \ref{niru-p2-t040}.
\end{pf}

\section{Partial sums of functions in $\mathcal{W}^0_{\mathcal{H}}(\alpha)$}

\begin{lem}\label{niru-p2-l020}

Let $f = h + \overline{g}$ be in the class $\mathcal{W}^0_{\mathcal{H}}(\alpha)$ with $\alpha \geq 0$. Then for $|\epsilon| = 1$ and $|z|< 1/2$,  we have
$$
{\rm Re\,}((s_3(h) + \epsilon s_3(g))' + \alpha z(s_3(h) + \epsilon s_3(g))'')> 1/4.
$$

\end{lem}

\begin{pf} Let  $f = h + \overline{g}\in\mathcal{W}^0_{\mathcal{H}}(\alpha)$. Then by Theorem \ref{niru-p2-t005}, $h + \epsilon g$ is in the class $\mathcal{W}(\alpha)$ for  $\epsilon~ (|\epsilon | =1)$.
Then ${\rm Re\,}F_{\epsilon}(z) > 0$ where
$$
F_{\epsilon}(z) = (h + \epsilon g)' + \alpha z (h + \epsilon g)'' = 1 +\sum_{k =1}^{\infty}c_k z^k.
$$
It is easy to see that
$$
|2c_2 - c_1^2|\leq 4 - |c_1|^2.
$$
Now let  $2 c_2 - c_1^2 = c$. Then $c_2 = {c}/{2} + {c_1^2}/{2}$ and $|c|\leq 4 - |c_1|^2.$
Let $c_1 z = \alpha + i \beta$ and $\surd{c} z = \gamma + i \delta,$ where $\alpha$, $\beta$, $\gamma$ and $\delta$ are real numbers. Then for $|z|< 1/2$ it is easy to see that
$$
\alpha^2 + \beta ^2 = |c_1|^2|z|^2 < \frac{|c_1|^2}{4}
$$
and
\begin{equation}\label{niru-p2-e040}
\gamma ^2 + \delta ^2 = |c||z|^2 < \frac{|c|}{4} \leq 1 - \frac{|c_1|^2}{4} \leq 1 - (\alpha^2 + \beta ^2).
\end{equation}

Therefore from (\ref{niru-p2-e040}) and $c_2 = {c}/{2} + {c_1^2}/{2}$, we obtain
\begin{align*}
&{\rm Re\,}((s_3(h) + \epsilon s_3(g))' + \alpha z(s_3(h) + \epsilon s_3(g))'')\\
& = {\rm Re\,} (1 + c_1 z + c_2 z^2)\\
& =  {\rm Re\,} (1 + c_1 z + \frac{c}{2} z^2 + \frac{c_1^2}{2} z^2)\\
& = 1 + \alpha + \left(\frac{\gamma^2}{2} - \frac{\delta ^2}{2}\right) + \left(\frac{\alpha^2}{2}- \frac{\beta^2}{2}\right)\\
& > 1 + \alpha + \left(\frac{\gamma^2}{2} - \frac{1 - (\alpha^2 + \beta^2 + \gamma^2)}{2}\right) +  \left(\frac{\alpha^2}{2}- \frac{\beta^2}{2}\right)\\
& = \frac{1}{4} + \left( \alpha + \frac{1}{2}\right)^2 + \gamma^2 \geq \frac{1}{4}
\end{align*}
which yields the desired result.
\end{pf}

\begin{thm}

Let $f\in \mathcal{W}^0_{\mathcal{H}}(\alpha)$. Then for each $q\geq 2$, $s_{1,q}(f)\in\mathcal{W}^0_{\mathcal{H}}(\alpha)$ for $|z|< 1/2.$

\end{thm}

\begin{pf}
Let $f  = h + \overline{g}$ be in the class $\mathcal{W}^0_{\mathcal{H}}(\alpha)$, where $h$ and $g$ are of the form (\ref{niru-p2-e001a}). Clearly,
$$
s_{1,q}(f)(z) = s_1(h)(z) + \overline{s_q(g)(z)} = z + \overline {\sum_{n=2}^{q} b_n z^n}.
$$
and ${\rm Re\,}(s_1(h)'(z) + \alpha z(s_1(h)''(z))) =1$. An application of Theorem \ref{niru-p2-t010} shows that
\begin{align*}
\noindent |s_q(g)'(z) + \alpha z s_q(g)''(z)| &= \left|\sum_{n =2}^{q}(n^2 \alpha + n (1 - \alpha)) b_n z ^{n-1}\right|\\
& \leq \sum_{n =2}^{q}(n^2 \alpha + n (1 - \alpha))| b_n| |z ^{n-1}| \leq \sum_{n =2}^{q}|z| ^{n-1}\\
& = \frac{|z|(1 - |z|^{q-1})}{1 - |z|} < \frac{|z|}{1 - |z|} < 1 \\
& =  {\rm Re\,}(s_1(h)'(z) + \alpha z(s_1(h)''(z)))
\end{align*}
when $|z|< 1/2.$
Therefore  $s_{1,q}(f)\in\mathcal{W}^0_{\mathcal{H}}(\alpha)$ in  $|z|< 1/2.$

\end{pf}

\begin{thm}
Let $f\in\mathcal{W}^0_{\mathcal{H}}(\alpha)$ and $p$ and $q$ satisfy one of the following conditions :
\begin{itemize}
\item[(i)] $3 \leq p < q$,
\item[(ii)] $p = q \geq 2$,
\item[(iii)] $p> q \geq 3$,
\item[(iv)] p= $3$ and q = $2$.
\end{itemize}
Then $s_{p,q}(f)\in\mathcal{W}^0_{\mathcal{H}}(\alpha)$ in  $|z|< 1/2$.

\end{thm}

\begin{pf}
Let $f  = h + \overline{g}$ be in the class $\mathcal{W}^0_{\mathcal{H}}(\alpha)$, where $h$ and $g$ are of the form (\ref{niru-p2-e001a}). So $\sigma_p(h)(z) = \sum_{k = p + 1}^{\infty}a_k z^k$ and $\sigma_q(g)(z) = \sum_{k = q + 1}^{\infty}b_k z^k$. Then
$h = s_p(h) + \sigma_p(h)$ and $g = s_q(g) + \sigma_q(g)$. To prove $s_{p,q}(f)\in\mathcal{W}^0_{\mathcal{H}}(\alpha)$ it sufficies to prove that
$ s_p(h) + \epsilon s_q(g)$ is in the class $\mathcal{W}(\alpha)$ for each $\epsilon~ (|\epsilon| = 1).$
If $f\in\mathcal{W}^0_{\mathcal{H}}(\alpha)$ then
\begin{align}\label{niru-p2-e045}
&{\rm Re\,}((s_p(h) + \epsilon s_q(g))' + \alpha z (s_p(h) + \epsilon s_q(g))'')\\\nonumber
&\geq{\rm Re\,}((h + \epsilon g)' + \alpha z ((h + \epsilon g)'') - |((\sigma_p(h) + \epsilon \sigma_q(g))' + \alpha z (\sigma_p(h) + \epsilon \sigma_q(g))'')|\\\nonumber
& \geq \frac{1 + |z|}{1 - |z|} - |((\sigma_p(h) + \epsilon \sigma_q(g))' + \alpha z (\sigma_p(h) + \epsilon \sigma_q(g))'')|
\end{align}

{\bf Case (i):}  $3 \leq p < q$

An application of Theorem \ref{niru-p2-t010} and Theorem \ref{niru-p2-t015} gives
\begin{align}\label{niru-p2-e055}
&|(\sigma_p(h) + \epsilon \sigma_q(g))' + \alpha z (\sigma_p(h) + \epsilon \sigma_q(g))''|\\\nonumber
& = \left|\frac{}{}((p + 1)+ \alpha p (p+1)) a_{p+1}z^p +...+ (q + \alpha q (q-1)) a_q z^{q-1}\right. \\\nonumber
&\quad\quad + \left.\sum_{k = q+1}^{\infty} (k + \alpha k (k-1))(a_k + \epsilon b_k)z^{k-1}\right| \\\nonumber
& \leq 2 (|z|^p + ... + |z|^{q -1} + 2 \sum_{k = q+1}^{\infty}|z|^{k - 1} ) = \frac{2 |z|^p}{1 - |z|}.
\end{align}
Using (\ref{niru-p2-e055}) in(\ref{niru-p2-e045}) we obtain
\begin{align}\label{niru-p2-e060}
&{\rm Re\,}((s_p(h) + \epsilon s_q(g))' + \alpha z (s_p(h) + \epsilon s_q(g))'')\\\nonumber
&\geq \frac{1 + |z|}{1 - |z|} - \frac{2 |z|^p}{1 - |z|}.
\end{align}

Now for $4 \leq p < q $ and $|z|= 1/2$, the inequality (\ref{niru-p2-e060}) yields
$$
{\rm Re\,}((s_p(h) + \epsilon s_q(g))' + \alpha z (s_p(h) + \epsilon s_q(g))'') \geq \frac{1}{3} - \frac{1}{4} > 0.
$$

Since ${\rm Re\,}((s_p(h) + \epsilon s_q(g))' + \alpha z (s_p(h) + \epsilon s_q(g))'')$ is harmonic, it assumes the minimum value on the circle $|z| = 1/2$. Therefore for  $4 \leq p < q,$ $s_{p,q}(f)\in\mathcal{W}^0_{\mathcal{H}}(\alpha)$ in  $|z|< 1/2$.

If $p =3 < q$ then an application of Theorem \ref{niru-p2-t010} and Lemma \ref{niru-p2-l020} shows that
\begin{align*}
&{\rm Re\,}((s_3(h) + \epsilon s_q(g))' + \alpha z (s_3(h) + \epsilon s_q(g))'')\\
& = {\rm Re\,}\left((s_3(h) + \epsilon s_3(g))' + \alpha z (s_3(h) + \epsilon s_3(g))'' + \epsilon \sum_{k =4}^{q}(k + \alpha k (k -1))b_k z^{k-1}\right)\\
& \geq \frac{1}{4}- \left|\epsilon \sum_{k =4}^{q}(k + \alpha k (k -1))b_k z^{k-1}\right|\\
& > \frac{1}{4} - \frac{|z|^3}{1 - |z|}.
\end{align*}
It is easy to see that
$$
{\rm Re\,}((s_3(h) + \epsilon s_q(g))' + \alpha z (s_3(h) + \epsilon s_q(g))'') > 0
$$
for $|z|< 1/2$ and hence $s_{3, q}(f)\in\mathcal{W}^0_{\mathcal{H}}(\alpha)$ for  $|z|< 1/2$.

{\bf Case (ii):} $p = q \geq 2$

If $ p = q \geq 4$, then the inequality (\ref{niru-p2-e060})  gives $s_{p, q}\in\mathcal{W}^0_{\mathcal{H}}(\alpha)$ in $|z|< 1/2$ and Lemma  \ref{niru-p2-l020} implies that $s_{3, 3}\in\mathcal{W}^0_{\mathcal{H}}(\alpha)$ in $|z|< 1/2$. For $ p = q = 2$, $s_{2, 2}(f)(z) = z + a_2 z^2 + \overline{b_2 z^2}$.  An application of Theorem \ref{niru-p2-t015} shows that
\begin{align*}
&{\rm Re\,}((s_2(h) + \epsilon s_2(g))' + \alpha z (s_2(h) + \epsilon s_2(g))'')\\
& = 1 + 2 (1 + \alpha){\rm Re\,}(a_2 + \epsilon b_2) z\\
& \geq 1 - 2 (1 + \alpha) |(a_2 + \epsilon b_2) z|\\
& \geq 1 - 2|z| > 0
\end{align*}
in the disk $|z|< 1/2$. Hence, $s_{2, 2}(f)\in\mathcal{W}^0_{\mathcal{H}}(\alpha)$ for  $|z|< 1/2$.

{\bf Case (iii):} $p> q \geq 3$

By Theorem \ref{niru-p2-t010} and Theorem \ref{niru-p2-t015}, we obtain
\begin{align}\label{niru-p2-e065}
&|(\sigma_p(h) + \epsilon \sigma_q(g))' + \alpha z (\sigma_p(h) + \epsilon \sigma_q(g))'')|\\\nonumber
& = \left|\frac{}{}(q + 1)+ \alpha q (q+1)) \epsilon b_{q+1}z^q +...+ (p + \alpha p (p-1)) \epsilon b_p z^{p-1} \right.\\\nonumber
& \quad\quad + \left.\sum_{k = p+1}^{\infty} (k + \alpha k (k-1))(a_k \epsilon b_k)z^{k-1}\right| \\\nonumber
& = \frac{|z|^q - |z|^p}{1 - |z|} + 2\frac{|z^p|}{1 - |z|} = \frac{|z|^q + |z|^p}{1 - |z|}.
\end{align}
Using (\ref{niru-p2-e065}) in (\ref{niru-p2-e045}) we obtain
\begin{align}\label{niru-p2-e070}
&{\rm Re\,}(s_p(h) + \epsilon s_q(g))' + \alpha z (s_p(h) + \epsilon s_q(g))'' \\\nonumber
& \geq \frac{1-|z|}{1 + |z|} - \frac{|z|^q + |z|^p}{1 - |z|}
\end{align}

If $q \geq 4$  then as in  Case (i) we obtain 
\begin{align*}
&{\rm Re\,}((s_p(h) + \epsilon s_q(g))' + \alpha z (s_p(h) + \epsilon s_q(g))'')\\\nonumber
&\geq \frac{1 + |z|}{1 - |z|} - \frac{2 |z|^q}{1 - |z|} > 0.
\end{align*}
Hence for  $q \geq 4$ and $p > q$, $s_{p, q}\in\mathcal{W}^0_{\mathcal{H}}(\alpha)$ for $|z|< 1/2$.

If $q = 3$ and $|z| = 1/2$,  from (\ref{niru-p2-e070}) it follows that
$$
{\rm Re\,}((s_p(h) + \epsilon s_3(g))' + \alpha z (s_p(h) + \epsilon s_3(g))'') \geq \frac{1}{12} - \frac{1}{2^{p-1}}.
$$
If $p \geq 5$ then the last estimate shows that
$$
{\rm Re\,}((s_p(h) + \epsilon s_3(g))' + \alpha z (s_p(h) + \epsilon s_3(g))'') >0
$$
and hence $s_{p, q}\in\mathcal{W}^0_{\mathcal{H}}(\alpha)$ for $|z|< 1/2$.

If $q = 3$ and $p = 4$, it follows that
\begin{align}\label{niru-p2-e075}
&{\rm Re\,}((s_4(h) + \epsilon s_3(g))' + \alpha z (s_4(h) + \epsilon s_3(g))'')\\\nonumber
&= {\rm Re\,}((s_4(h) + \epsilon s_3(g))' + \alpha z (s_4(h) + \epsilon s_3(g))'') + {\rm Re\,}((4 + 12 \alpha) a_4 z^3)\\\nonumber
&\geq  {\rm Re\,}((s_4(h) + \epsilon s_3(g))' + \alpha z (s_4(h) + \epsilon s_3(g))'') - |((4 + 12 \alpha) a_4 z^3)|.
\end{align}
Using Theorem \ref{niru-p2-t015} and Lemma \ref{niru-p2-l020} in (\ref{niru-p2-e075}) we obtain
\begin{align*}
&{\rm Re\,}((s_4(h) + \epsilon s_3(g))' + \alpha z (s_4(h) + \epsilon s_3(g))'')\\
&>\frac{1}{4} -2 |z|^3 > \frac{1}{4} -\frac{1}{2^2} = 0 \quad \mbox{for} \quad |z|= \frac{1}{2}.
\end{align*}
This implies $s_{4, 3}\in\mathcal{W}^0_{\mathcal{H}}(\alpha)$ for $|z|< 1/2$ and hence for $p> q \geq 3$, $s_{p, q}\in\mathcal{W}^0_{\mathcal{H}}(\alpha)$ for $|z|< 1/2$.

{\bf Case (iv):} p= $3$ and q = $2$

An application of Theorem \ref{niru-p2-t010} and Lemma \ref{niru-p2-l020} shows that
\begin{align*}
&{\rm Re\,}((s_3(h) + \epsilon s_2(g))' + \alpha z (s_3(h) + \epsilon s_2(g))'')\\
& ={\rm Re\,}((s_3(h) + \epsilon s_3(g))' + \alpha z (s_3(h) + \epsilon s_3(g))'') - {\rm Re\,}(((3 + 6 \alpha)\epsilon)b_3 z^2)\\
& > \frac{1}{4} - |z|^2 >  \frac{1}{4} - \frac{1}{2^2} = 0 \quad \mbox{for}\quad |z|= 1/2.
\end{align*}
Hence  $s_{3, 2}\in\mathcal{W}^0_{\mathcal{H}}(\alpha)$ for $|z|< 1/2.$ This completes the proof.
\end{pf}

\begin{thm}
Let $f\in\mathcal{W}^0_{\mathcal{H}}(\alpha)$. If $p = 2 < q,$ then $s_{2,q}(f)\in\mathcal{W}^0_{\mathcal{H}}(\alpha)$ in  $|z|< \frac{3 - \surd5}{2}.$
If $p\geq 4$ and $q =2$, then $s_{p,2}(f)\in\mathcal{W}^0_{\mathcal{H}}(\alpha)$ in  $|z|< r_{0}$ where $r_{0} \approx 0.433797$ is the unique real root of the equation $1 - 2r -r^3 - r^4 - r^5 = 0.$
\end{thm}

\begin{pf}
Let $p = 2 < q$. Then $s_{2, q}(f) = s_2 (h) + \overline{s_q(g)} = z + a_2 z^2 + \overline{\sum_{k = 2}^{q}b_n z^n}.$
So
\begin{align*}
&(s_2(h) + \epsilon s_q(g))' + \alpha z (s_2(h) + \epsilon s_q(g))''\\
& = 1 + 2a_2z(1 + \alpha) + \epsilon \sum_{k=2}^{q}(k + \alpha k(k-1))b_k z^{k-1}.
\end{align*}
An application of Theorem \ref{niru-p2-t015} shows that 
\begin{align*}
&|((s_2(h) + \epsilon s_q(g))' + \alpha z (s_2(h) + \epsilon s_q(g))'') - 1|\\
& = \left|2a_2z(1 + \alpha) + \epsilon \sum_{k=2}^{q}(k + \alpha k(k-1))b_k z^{k-1}\right|\\
& = \left|2(1 + \alpha(a_2 + \epsilon b_2)z + \epsilon \sum_{k=3}^{q}(k + \alpha k(k-1))b_k z^{k-1}\right|\\
& < 2 |z| + \frac{|z|^2}{1 - |z|}\\
& <1
\end{align*}
for  $|z|< \frac{3 - \surd5}{2}.$ Hence $s_{2,q}(f)\in\mathcal{W}^0_{\mathcal{H}}(\alpha)$ for  $|z|< \frac{3 - \surd5}{2}.$

Consider the case $p\geq 4$ and $q =2$. Then using (\ref{niru-p2-e070}) we obtain
\begin{align*}
&{\rm Re\,}((s_p(h) + \epsilon s_2(g))' + \alpha z (s_p(h) + \epsilon s_2(g))'')\\
&\geq \frac{1-|z|}{1 + |z|} - \frac{|z|^q + |z|^p}{1 - |z|}\\
&\geq \frac{1-|z|}{1 + |z|} - \frac{|z|^2 + |z|^4}{1 - |z|}\\
& = \frac{1 - 2|z| - |z|^3 - |z|^4 - |z|^5}{1 - |z|^2} > 0
\end{align*}
for $|z|< r_{0} \approx 433797$. Hence $s_{p,2}(f)\in\mathcal{W}^0_{\mathcal{H}}(\alpha)$ for $|z|< 0.433797$.

\end{pf}

\end{document}